\documentclass[11pt]{amsart}
\usepackage[margin=1.15in]{geometry}
\usepackage{amsmath,amssymb,amsthm,mathtools,enumitem}
\usepackage[colorlinks=true,linkcolor=blue,citecolor=blue,urlcolor=blue]{hyperref}

\usepackage{tikz-cd}

\numberwithin{equation}{section}
\newtheorem{theorem}{Theorem}[section]
\newtheorem{proposition}[theorem]{Proposition}
\newtheorem{lemma}[theorem]{Lemma}

\newtheorem*{maintheorem}{Main Theorem}
\newtheorem*{maincorollary}{Corollary}
\theoremstyle{definition}
\newtheorem{definition}[theorem]{Definition}

\theoremstyle{definition}
\newtheorem{remark}[theorem]{Remark}
\newtheorem{rex}[theorem]{Remark and Example}
\newtheorem{example}[theorem]{Example}

\newcommand{\R}{\mathbb{R}}
\newcommand{\Z}{\mathbb{Z}}
\newcommand{\lT}{\ell_T}
\newcommand{\leT}{\le_T}
\newcommand{\Red}{\operatorname{Red}}
\newcommand{\Comp}{\operatorname{Comp}}
\newcommand{\pr}{\operatorname{pr}}
\newcommand{\Hom}{\operatorname{Hom}}
\newcommand{\Ext}{\operatorname{Ext}}
\newcommand{\End}{\operatorname{End}}
\newcommand{\modd}{\operatorname{mod}}
\newcommand{\Wt}{\widetilde W}
\newcommand{\St}{\widetilde S}
\newcommand{\Tt}{\widetilde T}
\newcommand{\Vt}{\widetilde V}
\newcommand{\Bt}{\widetilde B}
\newcommand{\Phit}{\widetilde\Phi}
\newcommand{\at}{\widetilde\alpha}
\newcommand{\Nt}{\widetilde N}

\title[Prefixes of Coxeter elements]{Prefixes of Coxeter elements are Coxeter elements\\ of the reflection subgroups they generate}

\author[P.~Wegener]{Patrick Wegener}
\address{Patrick Wegener, Philipps-Universit\"at Marburg -- Studienkolleg Mittelhessen, Germany}
\email{wegenerp@staff.uni-marburg.de}

\date{\today}

\begin{document}

\begin{abstract}
Let $(W,S)$ be a Coxeter system of finite rank, let $T$ be its set of reflections and let $c \in W$ be a Coxeter element.
For $x$ in the interval $[1,c]_T$ of the absolute order, let $P(x)$ be the subgroup generated by the reflections $t\le_T x$.
We prove that $x$ is a Coxeter element of $P(x)$ with respect to the canonical Coxeter generators of $P(x)$, that $P(x)$ is generated by the factors of any reduced reflection factorization of $x$, and that the Hurwitz action on the reduced reflection factorizations of $x$ is transitive.
The proof reduces the general case to the case of a Coxeter system whose Coxeter matrix has entries in $\{1,2,3,\infty\}$, which was treated by Hubery and Krause.
\end{abstract}

\maketitle

\section{Introduction}\label{sec:intro}

Usually, I would now start to describe what this paper's main theorem is about and why it is important. I would explain that it might be an important part of the so-called \emph{dual approach to Coxeter and Artin groups} and that it could one day be used to prove the \emph{$K(\pi,1)$-conjecture for Artin groups} -- things many other mathematicians working in this area have written in slightly different phrasings before.

But this conjecture has now (most likely) been settled by OpenAI's models \cite{OpenAI} (although their proof does not use the dual approach). And to be honest, although parts of my mathematical career were dedicated to proving the main theorem of this paper, I was only able to prove it with massive AI assistance provided by OpenAI's ChatGPT Astra.

The $K(\pi,1)$-conjecture was so influential for an entire mathematical field and for my own academic career that, on the one hand, one might feel sorry that it was solved overnight. On the other hand, the dual approach to Coxeter and Artin groups evolved from the idea of proving the $K(\pi,1)$-conjecture, and it still exists regardless of whether the conjecture has been proven by other means.

So, in the end, I will not do what a mathematical paper is expected to start with: promoting the importance of its main theorem. Perhaps it has lost part of its significance due to the proof of the $K(\pi,1)$-conjecture. But, from my perspective, it remains a great piece of mathematics. Let me state it.

Let $(W,S)$ be a Coxeter system with $S$ finite, and let $T=\{wsw^{-1}\mid w\in W,\ s\in S\}$ be its set of reflections.
The \emph{reflection length} $\lT(w)$ of $w\in W$ is the smallest integer $k\ge 0$ such that $w$ is a product of $k$ reflections, and the \emph{absolute order} on $W$ is defined by
\[
u\leT v\quad:\Longleftrightarrow\quad \lT(u)+\lT(u^{-1}v)=\lT(v).
\]
A \emph{Coxeter element} of $(W,S)$ is a product of all elements of $S$ in some order. Note that such an element is often called a standard Coxeter element in the literature.  
For a Coxeter element $c$ we write $[1,c]_T=\{x\in W\mid x\leT c\}$.
For $x\in W$ with $\lT(x)=k$ let $\Red_T(x)$ be the set of $k$-tuples $(t_1,\dots,t_k)\in T^k$ with $x=t_1\cdots t_k$.
The braid group $\mathrm{Br}_k$ on $k$ strands acts on $\Red_T(x)$ by the \emph{Hurwitz action}
\begin{align*}
\sigma_i\cdot(t_1,\dots,t_k) & =(t_1,\dots,t_{i-1},\,t_{i+1},\,t_{i+1}t_it_{i+1},\,t_{i+2},\dots,t_k),\\
\sigma_i^{-1}\cdot(t_1,\dots,t_k) & =(t_1,\dots,t_{i-1},\,t_it_{i+1}t_i,\,t_i,\,t_{i+2},\dots,t_k) \qquad(1\le i<k).
\end{align*}
This also yields an action of the braid group on $k$-tuples of elements of $T$, and we use the term \emph{Hurwitz orbit} for its orbits.

By a theorem of Dyer and Deodhar \cite{De89,Dy90}, every reflection subgroup $W'$ of $W$ (that is, every subgroup generated by reflections) is a Coxeter group with respect to a canonical set $\chi(W')\subseteq T$ of Coxeter generators (see Theorem~\ref{thm:dyer} below). For $x\in W$ we put
\[
P(x)=\langle t\in T\mid t\leT x\rangle .
\]
We address the following two classical questions for a Coxeter system $(W,S)$, a Coxeter element $c \in W$ and $x\in[1,c]_T$ (see \cite[Questions~3.1 to~3.3]{Pa} and \cite[Questions~6.1 and~6.2]{G}).
\begin{itemize}
    \item Is $x$ a Coxeter element in $P(x)$?
    \item Is the Hurwitz action transitive on the $T$-reduced expressions of $x$?
\end{itemize}

As recalled in \cite[Remark~3.8]{G}, positive answers were previously known for finite Coxeter groups \cite[Section~1.4]{Be03}, for affine Coxeter groups \cite[Theorem~3.22]{PS21}, for Coxeter systems of rank three \cite[Theorem~3.10]{DPS}, and, by methods from the representation theory of hereditary algebras, for Weyl groups of symmetrisable Kac--Moody algebras, in particular for crystallographic Coxeter groups \cite[Theorem~5.2 and Corollary~5.8]{HK16}.
For $x=c$, Hurwitz transitivity on $\Red_T(c)$ is due to Igusa and Schiffler \cite[Theorem~1.4]{IS10}, and for Coxeter elements of parabolic subgroups it is due to Baumeister, Dyer, Stump and the author \cite[Theorem~1.3]{BDSW14}. Our main result gives an affirmative answer to both questions.

\begin{maintheorem}
Let $(W,S)$ be a Coxeter system of finite rank, $c\in W$ be a Coxeter element, and $x\in[1,c]_T$. Put $k=\lT(x)$ and let $(t_1,\dots,t_k)\in\Red_T(x)$.
Then the following hold.
\begin{enumerate}[label=\textup{(\roman*)}]
\item $P(x)=\langle t_1,\dots,t_k\rangle$. In particular, every reflection $t\leT x$ lies in $\langle t_1,\dots,t_k\rangle$.
\item $|\chi(P(x))|=k$, and $x=r_1\cdots r_k$, where $r_1,\dots,r_k$ are the $k$ elements of $\chi(P(x))$ in a suitable order. That is, $x$ is a Coxeter element of the Coxeter system $(P(x),\chi(P(x)))$.
\item The Hurwitz action on $\Red_T(x)$ is transitive.
\end{enumerate}
\end{maintheorem}

Since conjugation by $w \in W$ preserves the absolute length $\ell_T$, maps $\Red_T(x)$ onto $\Red_T(wxw^{-1})$ and commutes with the Hurwitz action, part (i) and (iii) of the Main Theorem hold more generally. 

\begin{maincorollary}
    Let $c$ be a Coxeter element of $(W,S)$ and $w \in W$. Then statements (i) and (iii) of the Main Theorem hold for every $x\in[1,wcw^{-1}]_T$ and every $(t_1,\dots,t_k)\in\Red_T(x)$.
\end{maincorollary}

The proof of the Main Theorem is organized as follows. Section~\ref{sec:prelim} fixes notation and recalls results that we will use later. Section~\ref{sec:subgroups} contains two elementary lemmas on reflection subgroups. Section~\ref{sec:crys} deals with Coxeter systems whose Coxeter matrix has all off-diagonal entries in $\{2,3,\infty\}$. For these, the transitivity of the Hurwitz action on $\Red_T(y)$ for $y$ below a Coxeter element is due to Hubery and Krause \cite{HK16}, and we deduce from their results that such a $y$ is a Coxeter element of the subgroup generated by the factors of any of its reduced reflection factorizations, with respect to the canonical Coxeter generators (Proposition~\ref{prop:crys}). In Section~\ref{sec:unfolding} we construct, for an arbitrary Coxeter system $(W,S)$ of finite rank, a homomorphism $F$ from $W$ into a Coxeter group $\Wt$ of this special type, which maps each reflection to a product of $\nu$ pairwise commuting reflections, for a fixed integer $\nu\ge1$.
We show that $F$ maps $[1,c]_T$ into $[1,F(c)]_{\Tt}$, multiplies the reflection length of the elements of $[1,c]_T$ by $\nu$, and satisfies the inequality $\nu\,\ell_{\chi(L)}(w)\le \ell_{\chi(H)}(F(w))$ for suitable pairs of reflection subgroups $L\subseteq W$ and $H\subseteq\Wt$ (Proposition~\ref{prop:transfer}).
In Section~\ref{sec:proof} we combine these results to prove the Main Theorem.

\medskip
\noindent \textbf{Disclosure of AI tools.} Generative AI tools, namely ChatGPT Astra and Claude Opus 5.5, were used during the exploratory and preparatory stages of this work. The author directed the mathematical strategy and proofs. AI assisted in several aspects like the formulation of auxiliary results and the elaboration of technical details, as well as consistency checks, and \LaTeX{} preparation. All AI-generated suggestions were independently reviewed and verified by the author, who assumes full responsibility for the mathematical arguments, results, and final content of the paper.

\section{Preliminaries}\label{sec:prelim}

Throughout, $(W,S)$ is a Coxeter system of finite rank with Coxeter matrix $(m_{st})_{s,t\in S}$, and $T$ is its set of reflections.
For a subset $X$ of a group we write $\langle X\rangle$ for the subgroup generated by $X$, and for a tuple $A=(a_1,\dots,a_k)$ of group elements we write $\langle A\rangle=\langle a_1,\dots,a_k\rangle$.

\subsection{Roots}
Let $V$ be the real vector space with basis $\{\alpha_s ~|~ s\in S\}$, equipped with the symmetric bilinear form $B$ given by 
$$B(\alpha_s,\alpha_t) = \begin{cases}
    -\cos(\pi/m_{st}), ~\text{if } m_{st}<\infty\\
    -1, ~\text{if } m_{st}=\infty
\end{cases}
$$
The group $W$ acts on $V$ by 
$$s(v)=v-2B(\alpha_s,v)\alpha_s$$ 
for $s\in S$, and $B$ is $W$-invariant. Let 
$$\Phi=\{w(\alpha_s)\mid w\in W,\ s\in S\}$$ 
be the set of roots, let $\Phi^+$ be the set of roots all of whose coordinates with respect to the basis $(\alpha_s)_{s\in S}$ are nonnegative, and put $\Phi^-=-\Phi^+$.

\begin{theorem}[{\cite[Sections~5.3--5.7]{Hum90}}]\label{thm:roots}
The following hold.
\begin{enumerate}[label=\textup{(\alph*)}]
\item The representation $W\to GL(V)$ is faithful.
\item $\Phi=\Phi^+\sqcup\Phi^-$.
\item If $\beta\in\Phi$ and $\lambda\in\R$ with $\lambda\beta\in\Phi$, then $\lambda=\pm1$.
\item For $\beta\in\Phi$ there is a unique $s_\beta\in T$ acting on $V$ by $v\mapsto v-2B(\beta,v)\beta$, and $\beta\mapsto s_\beta$ is a bijection from $\Phi^+$ onto $T$. For $t\in T$ we denote by $\alpha_t\in\Phi^+$ the root with $s_{\alpha_t}=t$. Moreover $s_{-\beta}=s_\beta$ and $ws_\beta w^{-1}=s_{w(\beta)}$ for all $w\in W$ and $\beta\in\Phi$.
\item For $w\in W$, the set
\[
N(w)=\{t\in T\mid \ell_S(tw)<\ell_S(w)\}=\{t\in T\mid w^{-1}(\alpha_t)\in\Phi^-\}
\]
has exactly $\ell_S(w)$ elements.
\end{enumerate}
\end{theorem}

\subsection{Reflection length and the Hurwitz action}
We note the following elementary facts.

\begin{lemma}\label{lem:elementary}
Let $w\in W$ with $\lT(w)=m$.
\begin{enumerate}[label=\textup{(\alph*)}]
\item The Hurwitz action preserves the product of a tuple and the subgroup generated by its entries. It maps $\Red_T(w)$ to itself.
\item If $(t_1,\dots,t_m)\in\Red_T(w)$ and $0\le j\le m$, then $(t_1,\dots,t_j)\in\Red_T(t_1\cdots t_j)$ and $t_1\cdots t_j\leT w$.
\item If $x\leT w$, $\lT(x)=k$ and $(t_1,\dots,t_k)\in\Red_T(x)$, then there are $t_{k+1},\dots,t_m\in T$ with $(t_1,\dots,t_m)\in\Red_T(w)$.
\item If $(t_1,\dots,t_m)\in\Red_T(w)$, then $t_i\leT w$ for $1\le i\le m$. If $t\in T$ satisfies $t\leT w$, then there is an element of $\Red_T(w)$ whose first entry is $t$.
\end{enumerate}
\end{lemma}

\begin{lemma}[{\cite[Lemma~2.1]{BDSW14}}]\label{thm:lengthc}
If $c$ is a Coxeter element of $(W,S)$, then $\lT(c)=|S|$.
\end{lemma}

\subsection{Reflection subgroups}
A \emph{reflection subgroup} of $W$ is a subgroup generated by a subset of $T$.
For a reflection subgroup $W'$ of $W$ put
\[
\chi(W')=\{t\in T\mid N(t)\cap W'=\{t\}\}.
\]

\begin{theorem}[{\cite[Theorem~3.3 and Corollary~3.11]{Dy90}}, \cite{De89}]\label{thm:dyer}
Let $W'$ be a reflection subgroup of $W$.
\begin{enumerate}[label=\textup{(\alph*)}]
\item $(W',\chi(W'))$ is a Coxeter system, and $W'\cap T=\bigcup_{w\in W'}w\chi(W')w^{-1}$ is its set of reflections.
\item For every $w\in W'$ we have $\ell_{\chi(W')}(w)=|N(w)\cap W'|$, where $\ell_{\chi(W')}$ is the length function of $(W',\chi(W'))$.
\item If $W'=\langle R\rangle$ for a subset $R\subseteq T$, then $|\chi(W')|\le|R|$ and $W'\cap T=\bigcup_{w\in W'}wRw^{-1}$.
\end{enumerate}
\end{theorem}

The elements of $\chi(W')$ are called the \emph{canonical generators} of $W'$, and $|\chi(W')|$ is the \emph{rank} of $W'$.
We write $T(W')=W'\cap T$.
Note that $\chi(W')\subseteq T(W')\subseteq T$, so that
\begin{equation}\label{eq:lengthineq}
\lT(w)\le\ell_{T(W')}(w)\le\ell_{\chi(W')}(w)\qquad(w\in W'),
\end{equation}
where $\ell_{T(W')}$ denotes the reflection length of the Coxeter system $(W',\chi(W'))$.

\begin{theorem}[{\cite[Theorem~4.4]{Dy90}}]\label{thm:dyercrit}
Let $\Pi\subseteq\Phi^+$ be a finite set such that
\[
B(\beta,\gamma)\in\{-\cos(\pi/m)\mid m\in\Z,\ m\ge2\}\cup(-\infty,-1]\qquad\text{for all }\beta\ne\gamma\text{ in }\Pi .
\]
Then $\chi(\langle s_\beta\mid\beta\in\Pi\rangle)=\{s_\beta\mid\beta\in\Pi\}$.
\end{theorem}

\begin{theorem}[{\cite[Theorem~1.3]{BDSW14}}]\label{thm:bdsw}
Let $(W,S)$ be a Coxeter system of finite rank with set of reflections $T$, and let $w=s_1\cdots s_m$ with pairwise distinct $s_1,\dots,s_m\in S$. Then the Hurwitz action on $\Red_{T}(w)$ is transitive.
\end{theorem}

\section{Two lemmas on reflection subgroups}\label{sec:subgroups}

The following result follows immediately as a consequence of the Hurwitz action.
\begin{lemma}\label{lem:distinct}
Let $w\in W$ and $r_1,\dots,r_k\in T$ with $w=r_1\cdots r_k$ and $k=\lT(w)$. Then $r_i\ne r_j$ for all $1\le i<j\le k$.
\end{lemma}

\begin{lemma}\label{lem:parabolic}
Let $W'$ be a reflection subgroup of $W$ such that $\chi(W')$ is finite, and let $w\in W'$ with $\ell_{\chi(W')}(w)=\lT(w)=:k$.
\begin{enumerate}[label=\textup{(\alph*)}]
\item There are pairwise distinct $r_1,\dots,r_k\in\chi(W')$ with $w=r_1\cdots r_k$.
\item Let $r_1,\dots,r_k$ be as in \textup{(a)}. Every $A\in\Red_T(w)$ all of whose entries lie in $W'$ lies in the Hurwitz orbit of $(r_1,\dots,r_k)$. In particular $\langle A\rangle=\langle r_1,\dots,r_k\rangle$.
\end{enumerate}
\end{lemma}

\begin{proof}
(a) Let $w=r_1\cdots r_k$ be a reduced expression of $w$ in the Coxeter system $(W',\chi(W'))$. Since $\chi(W')\subseteq T$ and $k=\lT(w)$, the $r_i$ are pairwise distinct by Lemma~\ref{lem:distinct}.

(b) By \eqref{eq:lengthineq} we have $k=\lT(w)\le\ell_{T(W')}(w)\le\ell_{\chi(W')}(w)=k$, so $\ell_{T(W')}(w)=k$.
The entries of $A$ lie in $W'\cap T=T(W')$, and the $r_i$ lie in $\chi(W')\subseteq T(W')$; hence both $A$ and $(r_1,\dots,r_k)$ belong to $\Red_{T(W')}(w)$.
By Theorem~\ref{thm:dyer}(a), $(W',\chi(W'))$ is a Coxeter system of finite rank with set of reflections $T(W')$.
By (a) and Theorem~\ref{thm:bdsw}, applied to $(W',\chi(W'))$, the tuples $A$ and $(r_1,\dots,r_k)$ lie in the same orbit of the Hurwitz action on $\Red_{T(W')}(w)$,
hence in the same Hurwitz orbit in $\Red_T(w)$.
The last assertion follows from Lemma~\ref{lem:elementary}(a).
\end{proof}

\section{Coxeter systems with labels in \texorpdfstring{$\{2,3,\infty\}$}{\{2,3,infinity\}}}\label{sec:crys}

We call a Coxeter system \emph{integral} if all off-diagonal entries of its Coxeter matrix lie in $\{2,3,\infty\}$.
Throughout this section $(\Wt,\St)$ denotes an integral Coxeter system of finite rank $n$, with geometric representation $\Vt$, simple roots $\at_s$ $(s\in\St)$, bilinear form $\Bt$, roots $\Phit=\Phit^+\sqcup\Phit^-$ and set of reflections $\Tt$.
Note that $2\Bt(\at_s,\at_t)\in\{2,0,-1,-2\}$ for all $s,t\in\St$.

\begin{lemma}\label{lem:integral}
The lattice $\widetilde\Gamma=\bigoplus_{s\in\St}\Z\at_s$ is $\Wt$-stable. Every root $\beta\in\Phit$ lies in $\widetilde\Gamma$ and satisfies $\Bt(\beta,\beta)=1$.
\end{lemma}

\begin{proof}
For $s,t\in\St$ we have $s(\at_t)=\at_t-2\Bt(\at_s,\at_t)\at_s\in\widetilde\Gamma$, so $\widetilde\Gamma$ is $\Wt$-stable and contains $\Phit=\Wt\{\at_s\mid s\in\St\}$. The second assertion follows from $\Bt(\at_s,\at_s)=1$ and the $\Wt$-invariance of $\Bt$.
\end{proof}

\subsection{Generalised Cartan lattices}
We recall the following notions from \cite[Sections~2 and~3]{HK16}.
A \emph{bilinear lattice} is an abelian group $\Gamma$ together with a non-degenerate bilinear form $\langle-,-\rangle\colon\Gamma\times\Gamma\to\Z$ (all bilinear lattices occurring below are free of finite rank). Its symmetrisation is 
$$
(x,y)=\langle x,y\rangle+\langle y,x\rangle \quad \text{for }x,y \in \Gamma.
$$
An element $a\in\Gamma$ is a \emph{root} if $\langle a,a\rangle>0$ and $\langle a,x\rangle/\langle a,a\rangle,\ \langle x,a\rangle/\langle a,a\rangle\in\Z$ for all $x\in\Gamma$. Corresponding to a root $a$ we have the reflection 
$$s_a: \Gamma \rightarrow \Gamma, ~s_a(x)=x-2\frac{(x,a)}{(a,a)}a.
$$

A sequence $E=(e_1,\dots,e_r)$ of roots is \emph{exceptional} if $\langle e_i,e_j\rangle=0$ for all $i>j$; it is \emph{complete} if $e_1,\dots,e_r$ generate $\Gamma$, and \emph{orthogonal} if $\langle e_i,e_j\rangle\le0$ for all $i\ne j$.
A \emph{generalised Cartan lattice} is a pair $(\Gamma,E)$ consisting of a bilinear lattice $\Gamma$ and a complete orthogonal exceptional sequence $E=(e_1,\dots,e_n)$.
Its \emph{Weyl group} $W(\Gamma,E)$ is the subgroup of $\operatorname{Aut}(\Gamma)$ generated by $s_{e_1},\dots,s_{e_n}$, and $(W(\Gamma,E),\{s_{e_1},\dots,s_{e_n}\})$ is a Coxeter system \cite[Proposition~3.13]{Ka90}.
The set of \emph{real roots} is 
$$\Phi(\Gamma,E)=\{w(e_i)\mid w\in W(\Gamma,E),\ 1\le i\le n\},
$$
the reflections of $W(\Gamma,E)$ are the elements $s_a$ with $a\in\Phi(\Gamma,E)$, the Coxeter element $\mathrm{cox}(\Gamma)=s_{e_1}\cdots s_{e_n}$ does not depend on $E$ \cite[Proposition~2.4]{HK16}, and $\mathrm{NC}(\Gamma,E)$ denotes the interval $[1,\mathrm{cox}(\Gamma)]$ in the absolute order of $W(\Gamma,E)$.

\begin{lemma}\label{lem:gcl}
Let $\St=\{s_1,\dots,s_n\}$ and $d=s_1s_2\cdots s_n$. Let $\Gamma=\widetilde\Gamma$, $e_i=\at_{s_i}$ and $E=(e_1,\dots,e_n)$, and define a bilinear form on $\Gamma$ by
\[
\langle e_i,e_j\rangle=\begin{cases}2\Bt(e_i,e_j)&\text{if }i<j,\\ 1&\text{if }i=j,\\ 0&\text{if }i>j.\end{cases}
\]
Then $(\Gamma,E)$ is a generalised Cartan lattice with $(x,y)=2\Bt(x,y)$ for $x,y\in\Gamma$.
Restriction to $\Gamma$ defines an isomorphism $\Wt\to W(\Gamma,E)$ which maps $s_i$ to $s_{e_i}$, $d$ to $\mathrm{cox}(\Gamma)$, and every $t\in\Tt$ to $s_{\at_t}$.
Moreover $\Phi(\Gamma,E)=\Phit$, and under this isomorphism $[1,d]_{\Tt}$ corresponds to $\mathrm{NC}(\Gamma,E)$.
\end{lemma}

\begin{proof}
The Gram matrix of $\langle-,-\rangle$ with respect to $E$ is unitriangular, so the form is non-degenerate, and $(e_i,e_j)=2\Bt(e_i,e_j)$ for all $i,j$, hence $(x,y)=2\Bt(x,y)$.
Each $e_i$ is a root, since $\langle e_i,e_i\rangle=1$. By definition $E$ is exceptional, orthogonal (as $2\Bt(e_i,e_j)\le0$ for $i\ne j$ and the Coxeter system is integral) and complete.
By Lemma~\ref{lem:integral}, every $\beta\in\Phit$ lies in $\Gamma$ and satisfies $\langle\beta,\beta\rangle=\tfrac12(\beta,\beta)=\Bt(\beta,\beta)=1$. Hence $\beta$ is a root of $\Gamma$, and 
$$
s_\beta(x)=x-(x,\beta)\beta=x-2\Bt(\beta,x)\beta
$$
for $x\in\Gamma$.
By Theorem~\ref{thm:roots}(d), applied to $(\Wt,\St)$, the reflection $s_\beta$ of $\Gamma$ is therefore the restriction to $\Gamma$ of the reflection $s_\beta\in\Tt$.
Since $\Gamma$ is $\Wt$-stable, restriction to $\Gamma$ is a homomorphism $\Wt\to\operatorname{Aut}(\Gamma)$, which maps $s_i$ to $s_{e_i}$, so its image is $W(\Gamma,E)$.
It is injective because $\Gamma$ spans $\Vt$ and $\Vt$ is a faithful representation (Theorem~\ref{thm:roots}(a)).
It follows that $\Phi(\Gamma,E)=\Phit$, that $d$ is mapped to $\mathrm{cox}(\Gamma)$, and that every $t\in\Tt$ is mapped to $s_{\at_t}$.
Since the isomorphism $\Wt\to W(\Gamma,E)$ maps $\Tt$ onto the set of reflections of $W(\Gamma,E)$, it preserves absolute length and absolute order, so $[1,d]_{\Tt}$ corresponds to $\mathrm{NC}(\Gamma,E)$.
\end{proof}

\begin{theorem}[Hubery--Krause \cite{HK16}]\label{thm:hk}
Let $d$ be a Coxeter element of $(\Wt,\St)$ and let $y\in[1,d]_{\Tt}$. Then the Hurwitz action on $\Red_{\Tt}(y)$ is transitive.
\end{theorem}

Note that this result is stated explicitly in the introduction of \cite{HK16}. For the sake of completeness, we provide a proof here. 

\begin{proof}
Call a sequence of real roots a \emph{real exceptional sequence} if it is a subsequence of a complete exceptional sequence consisting of real roots.
Every element of $\Red_{\Tt}(y)$ is of the form $(s_{f_1},\dots,s_{f_q})$ for a real exceptional sequence $(f_1,\dots,f_q)$.
To see this, let $(g_1,\dots,g_q)\in\Red_{\Tt}(y)$. By Lemma~\ref{thm:lengthc} and Lemma~\ref{lem:elementary}(c) it extends to some $(g_1,\dots,g_n)\in\Red_{\Tt}(d)$, where $n$ is the rank of $(\Wt, \St)$, and by \cite[Theorem~1.4]{IS10} (see also Theorem~\ref{thm:bdsw}) there is a braid $\sigma\in\mathrm{Br}_n$ with $(g_1,\dots,g_n)=\sigma\cdot(s_{e_1},\dots,s_{e_n})$.
The sequence $(f_1,\dots,f_n)=\sigma\cdot (e_1,\dots, e_n)$ is a complete exceptional sequence consisting of real roots \cite[Proposition~2.6]{HK16}, so $(f_1,\dots,f_q)$ is a real exceptional sequence, and $s_{f_i}=g_i$ for all $i$ by \cite[Lemma~3.4]{HK16}.
By \cite[Theorem~5.2]{HK16}, two real exceptional sequences $F=(f_1,\dots,f_q)$ and $F'=(f'_1,\dots,f'_q)$ with $s_{f_1}\cdots s_{f_q}=s_{f'_1}\cdots s_{f'_q}$ lie in the same orbit of the group $\{\pm1\}^q\rtimes\mathrm{Br}_q$ acting on real exceptional sequences.
Finally, the map $F\mapsto(s_{f_1},\dots,s_{f_q})$ is $\mathrm{Br}_q$-equivariant \cite[Lemma~3.4]{HK16}, and it is invariant under sign changes since $s_{-f}=s_f$.
\end{proof}

\subsection{Hereditary algebras}
Let $\Bbbk$ be a field and $\Lambda$ a finite-dimensional hereditary $\Bbbk$-algebra, and let $\modd\Lambda$ be the category of finite-dimensional $\Lambda$-modules.
All subcategories are assumed to be full and closed under isomorphisms.
The Grothendieck group $K_0(\Lambda)$ is free with basis the classes of the simple $\Lambda$-modules, and by \cite[Lemma~4.1(1)]{HK16} it is a bilinear lattice with respect to the Euler form
\[
\langle [X],[Y]\rangle=\dim_\Bbbk\Hom_\Lambda(X,Y)-\dim_\Bbbk\Ext^1_\Lambda(X,Y).
\]

For $M\in\modd\Lambda$ we have $[M]=\sum_S[M:S]\,[S]$, where $S$ runs through a set of representatives of the simple $\Lambda$-modules and $[M:S]$ is the multiplicity of $S$ as a composition factor of $M$ \cite[Chapter~III, Theorem~3.5]{ASS06}.
A module $X$ is \emph{exceptional} if it is indecomposable and $\Ext^1_\Lambda(X,X)=0$.
A sequence $X=(X_1,\dots,X_r)$ of modules is \emph{exceptional} if each $X_i$ is exceptional and $\Hom_\Lambda(X_i,X_j)=0=\Ext^1_\Lambda(X_i,X_j)$ for all $i>j$.
For such $X$ let $\mathcal C(X)$ be the smallest subcategory of $\modd\Lambda$ containing $X_1,\dots,X_r$ and closed under kernels, cokernels and extensions.
The braid group $\mathrm{Br}_r$ acts on the set of isomorphism classes of exceptional sequences of length $r$ in $\modd\Lambda$ \cite{CB93, Ringel}.

\begin{theorem}\label{thm:hkmod}
Let $\Bbbk$ be a field and $\Lambda$ a finite-dimensional hereditary $\Bbbk$-algebra.
\begin{enumerate}[label=\textup{(\alph*)}]
\item The simple $\Lambda$-modules can be ordered as $S_1,\dots,S_n$ such that $([S_1],\dots,[S_n])$ is a complete orthogonal exceptional sequence in $K_0(\Lambda)$ \cite[Lemma~4.1(1)]{HK16}.
\item Let $(\Gamma,E)$ be a generalised Cartan lattice, $E=(e_1,\dots,e_n)$, and let $\Bbbk$ be finite. Then there exist a finite-dimensional hereditary $\Bbbk$-algebra $\Lambda$ and an isometry $\Gamma\xrightarrow{\ \sim\ }K_0(\Lambda)$ sending each $e_i$ to the class of a simple $\Lambda$-module \cite[Lemma~4.1(2)]{HK16}.
\end{enumerate}
In the situation of \textup{(a)}, write $W(K_0(\Lambda))$, $\Phi(K_0(\Lambda))$ and $\mathrm{NC}(K_0(\Lambda))$ for the Weyl group, the real roots and the non-crossing partitions of the generalised Cartan lattice $(K_0(\Lambda),([S_1],\dots,[S_n]))$. They do not depend on the ordering in \textup{(a)}, since $W(K_0(\Lambda))$ is generated by the reflections $s_{[S]}$, $S$ simple, and $\mathrm{cox}(K_0(\Lambda))$ does not depend on the choice of a complete exceptional sequence \cite[Proposition~2.4]{HK16}.
\begin{enumerate}[label=\textup{(\alph*)},resume]
\item If $X$ is an exceptional $\Lambda$-module, then $[X]$ is a real root \cite[Proposition~4.4]{HK16}.
\item For $r\ge1$, the map $(X_1,\dots,X_r)\mapsto(s_{[X_1]},\dots,s_{[X_r]})$ induces a $\mathrm{Br}_r$-equivariant bijection from the set of isomorphism classes of exceptional sequences of length $r$ in $\modd\Lambda$ onto the set of $r$-tuples of reflections of $W(K_0(\Lambda))$ whose product lies in $\mathrm{NC}(K_0(\Lambda))$ and has absolute length $r$ \cite[Proposition~4.6]{HK16}.
\item If $X$ is an exceptional sequence of length $r$ in $\modd\Lambda$, then there are a finite-dimensional hereditary algebra $\Lambda'$ over a field $\Bbbk'$ with $\operatorname{rk}K_0(\Lambda')=r$ and an equivalence of categories $\modd\Lambda'\xrightarrow{\ \sim\ }\mathcal C(X)$ \cite[Theorem~A.4 and p.~2290]{HK16}.
\item Two exceptional sequences $X$ and $Y$ in $\modd\Lambda$ lie in the same orbit of the braid group action if and only if $\mathrm{cox}(\mathcal C(X))=\mathrm{cox}(\mathcal C(Y))$, where $\mathcal C\mapsto\mathrm{cox}(\mathcal C)$ is the poset isomorphism from the subcategories of the form $\mathcal C(X)$ onto $\mathrm{NC}(K_0(\Lambda))$ of \cite[Corollary~7.5]{HK16}. In particular, if $X$ and $Y$ have the same length and $\mathcal C(X)=\mathcal C(Y)$, then $X$ and $Y$ lie in the same orbit of the braid group action.
\end{enumerate}
\end{theorem}

We recall two standard facts which we use for a subcategory $\mathcal C$ of $\modd\Lambda$ closed under kernels, cokernels and extensions.
First, $\mathcal C$ is an abelian category and the inclusion $\mathcal C\to\modd\Lambda$ is exact.
Second, $\Hom_{\mathcal C}(M,M')=\Hom_\Lambda(M,M')$ and $\Ext^1_{\mathcal C}(M,M')=\Ext^1_\Lambda(M,M')$ for $M,M'\in\mathcal C$, since $\mathcal C$ is full and closed under extensions (see \cite[proof of Lemma~7.1]{HK16}).

\begin{proposition}\label{prop:crys}
Let $d$ be a Coxeter element of $(\Wt,\St)$, let $y\in[1,d]_{\Tt}$, put $q=\ell_{\Tt}(y)$, let $(g_1,\dots,g_q)\in\Red_{\Tt}(y)$ and $H=\langle g_1,\dots,g_q\rangle$.
Then $|\chi(H)|=q$, and $y=h_1\cdots h_q$, where $h_1,\dots,h_q$ are the $q$ elements of $\chi(H)$ in a suitable order.
In particular $\ell_{\chi(H)}(y)=q$.
\end{proposition}

\begin{proof}
If $q=0$, then $y=1$, $H=\{1\}$ and $\chi(H)=\emptyset$, so there is nothing to prove. Let $q\ge1$.
Write $\St=\{s_1,\dots,s_n\}$ with $d=s_1\cdots s_n$, and let $(\Gamma,E)$ be the generalised Cartan lattice of Lemma~\ref{lem:gcl}.
Fix a finite field $\Bbbk$.
By Theorem~\ref{thm:hkmod}(b) there are a finite-dimensional hereditary $\Bbbk$-algebra $\Lambda$ and an isometry $\iota\colon\Gamma\to K_0(\Lambda)$ such that $\iota(e_i)=[S_i]$ for simple $\Lambda$-modules $S_i$.
Since $\iota$ is bijective and the classes of the simple modules form a basis of $K_0(\Lambda)$, the modules $S_1,\dots,S_n$ form a set of representatives of the isomorphism classes of simple $\Lambda$-modules \cite[Chapter~III, Theorem~3.5]{ASS06}, and $([S_1],\dots,[S_n])=\iota(E)$ is a complete orthogonal exceptional sequence in $K_0(\Lambda)$.

Since $\iota$ preserves the symmetrised forms, we have $\iota s_{\alpha}\iota^{-1}=s_{\iota(\alpha)}$ for every root $\alpha$ of $\Gamma$. Consequently $\iota$ transports $W(\Gamma,E)$, $\Phi(\Gamma,E)$, $\mathrm{cox}(\Gamma)$ and $\mathrm{NC}(\Gamma,E)$ to $W(K_0(\Lambda))$, the real roots of $K_0(\Lambda)$, $\mathrm{cox}(K_0(\Lambda))$ and $\mathrm{NC}(K_0(\Lambda))$, respectively.
Combining this with Lemma~\ref{lem:gcl}, we identify $\Vt\supseteq\Gamma=K_0(\Lambda)$ via $\iota$, so that $\at_{s_i}=[S_i]$ and $2\Bt(a,b)=(a,b)$ for $a,b\in K_0(\Lambda)$, and we identify $\Wt=W(K_0(\Lambda))$, $\Phit=\Phi(K_0(\Lambda))$ and $[1,d]_{\Tt}=\mathrm{NC}(K_0(\Lambda))$, where $t\in\Tt$ is identified with $s_{\at_t}$.
For every $M\in\modd\Lambda$ the class $[M]=\sum_i[M:S_i]\,\at_{s_i}$ has nonnegative coordinates, and it is nonzero if $M\ne0$.

By Theorem~\ref{thm:hkmod}(d) there is an exceptional sequence $X=(X_1,\dots,X_q)$ in $\modd\Lambda$ with $s_{[X_i]}=g_i$ for $1\le i\le q$.
By Theorem~\ref{thm:hkmod}(e) there are a finite-dimensional hereditary algebra $\Lambda'$ over a field $\Bbbk'$ with $\operatorname{rk}K_0(\Lambda')=q$ and an equivalence of categories $\Psi\colon\modd\Lambda'\to\mathcal C(X)$.
By Theorem~\ref{thm:hkmod}(a), applied to $\Lambda'$, the simple $\Lambda'$-modules can be ordered as $S'_1,\dots,S'_q$ such that $([S'_1],\dots,[S'_q])$ is an exceptional sequence in $K_0(\Lambda')$.
Put $M_j=\Psi(S'_j)$ and $\gamma_j=[M_j]\in K_0(\Lambda)$ for $1\le j\le q$.
Then $M_1,\dots,M_q$ represent the isomorphism classes of simple objects of the abelian category $\mathcal C(X)$.

We claim that $M=(M_1,\dots,M_q)$ is an exceptional sequence in $\modd\Lambda$, and $\Hom_\Lambda(M_i,M_j)=0$ for $i\ne j$.

Since $\Psi$ is an equivalence of abelian categories, it is exact and induces isomorphisms of abelian groups $\Hom_{\Lambda'}(S'_i,S'_j)\cong\Hom_{\mathcal C(X)}(M_i,M_j)=\Hom_\Lambda(M_i,M_j)$ and $\Ext^1_{\Lambda'}(S'_i,S'_j)\cong\Ext^1_{\mathcal C(X)}(M_i,M_j)=\Ext^1_\Lambda(M_i,M_j)$, for all $i,j$.
By Schur's lemma every nonzero $\Lambda'$-homomorphism between $S_i'$ and $S_j'$ must be an isomorphism. Hence we have $\Hom_{\Lambda'}(S'_i,S'_j)=0$ for $i \ne j$. In particular $\Hom_\Lambda(M_i,M_j)=0$, and for $i>j$ the Euler form of $K_0(\Lambda')$ gives
\[
\dim_{\Bbbk'}\Ext^1_{\Lambda'}(S'_i,S'_j)=-\langle[S'_i],[S'_j]\rangle=0,
\]
so that $\Ext^1_\Lambda(M_i,M_j)=0$.

Let $1\le i\le q$. The ring $D=\End_{\Lambda'}(S'_i)$ is a division algebra, and $\End_\Lambda(M_i)\cong D$. In particular $M_i$ is indecomposable. Since $\Ext^1_{\Lambda'}(S'_i,S'_i)$ is a left $D$-vector space, $\dim_{\Bbbk'}\Ext^1_{\Lambda'}(S'_i,S'_i)=\mu\dim_{\Bbbk'}D$ for some integer $\mu\ge0$.
As $[S'_i]$ is a root of $K_0(\Lambda')$, we get 
\begin{align*}
0 < \langle [S'_i],[S'_i]\rangle& =\dim_{\Bbbk'}\Hom_{\Lambda'}(S_i',S_i')-\dim_{\Bbbk'}\Ext^1_{\Lambda'}(S_i',S_i')\\
& = \dim_{\Bbbk'}D - \mu\dim_{\Bbbk'}D = (1-\mu)\dim_{\Bbbk'}D,
\end{align*}
so $\mu=0$ and $\Ext^1_\Lambda(M_i,M_i)\cong\Ext^1_{\Lambda'}(S'_i,S'_i)=0$. Thus each $M_i$ is exceptional, and $M$ is an exceptional sequence. Furthermore we obtain $\mathcal C(M)=\mathcal C(X)$ (see also \cite[proof of Theorem~A.4, (2)$\Rightarrow$(3)]{HK16}).

Hence we can apply Theorem~\ref{thm:hkmod}(f) to find a braid $\beta\in \mathrm{Br}_q$ such that $M$ is isomorphic to $\beta\cdot X$. By the $\mathrm{Br}_q$-equivariance in Theorem~\ref{thm:hkmod}(d),
\[
(s_{\gamma_1},\dots,s_{\gamma_q})=\beta\cdot(s_{[X_1]},\dots,s_{[X_q]})=\beta\cdot(g_1,\dots,g_q).
\]
By Lemma~\ref{lem:elementary}(a) we obtain $(s_{\gamma_1},\dots,s_{\gamma_q})\in\Red_{\Tt}(y)$, $y=s_{\gamma_1}\cdots s_{\gamma_q}$ and $\langle s_{\gamma_1},\dots,s_{\gamma_q}\rangle=H$.

Since $M_j$ is exceptional, Theorem~\ref{thm:hkmod}(c) yields that each $\gamma_j$ is a real root, that is, $\gamma_j\in\Phit$. We have shown above it has nonnegative coordinates and is nonzero, so $\gamma_j\in\Phit^+$.
As $(s_{\gamma_1},\dots,s_{\gamma_q})\in\Red_{\Tt}(y)$, Lemma~\ref{lem:distinct} yields that the reflections $s_{\gamma_1},\dots,s_{\gamma_q}$ are pairwise distinct, hence so are $\gamma_1,\dots,\gamma_q$.
For $i\ne j$, using $\Hom_{\Lambda}(M_i,M_j) = 0$ and evaluating the Euler form gives
\begin{align*}
2\Bt(\gamma_i,\gamma_j) \stackrel{(\text{Lemma~\ref{lem:gcl}})}{=} (\gamma_i, \gamma_j) & =\langle\gamma_i,\gamma_j\rangle+\langle\gamma_j,\gamma_i\rangle\\
& =-\dim_\Bbbk\Ext^1_\Lambda(M_i,M_j)-\dim_\Bbbk\Ext^1_\Lambda(M_j,M_i)\in\Z_{\le0},
\end{align*}
so that $\Bt(\gamma_i,\gamma_j)\in\{0,-\tfrac12\}\cup(-\infty,-1]=\{-\cos(\pi/2),-\cos(\pi/3)\}\cup(-\infty,-1]$.
By Theorem~\ref{thm:dyercrit},
\[
\chi(H)=\chi(\langle s_{\gamma_1},\dots,s_{\gamma_q}\rangle)=\{s_{\gamma_1},\dots,s_{\gamma_q}\},
\]
and $y=s_{\gamma_1}\cdots s_{\gamma_q}$. Finally, $q=\ell_{\Tt}(y)\le\ell_{\chi(H)}(y)\le q$ by \eqref{eq:lengthineq}.
\end{proof}

\section{Unfolding}\label{sec:unfolding}

In this section $(W,S)$ is an arbitrary Coxeter system of finite rank.
We construct a homomorphism $F$ from $W$ to an integral Coxeter group $\Wt$ which maps each simple reflection to the longest element of a finite standard parabolic subgroup of $\Wt$. Homomorphisms of this kind have been studied by M\"uhlherr \cite{Mu93} and Crisp \cite{Cr99}.

\subsection{Construction}
Let $\mathcal M=\{m_{st}\mid s,t\in S,\ s\ne t,\ 3\le m_{st}<\infty\}$ be the set of finite Coxeter labels of $(W,S)$ which are at least $3$, and put
\begin{align*}
I&=\prod_{m\in\mathcal M}\{1,\dots,m-1\}=\bigl\{a\colon\mathcal M\to\Z\ \big|\ 1\le a(m)\le m-1\text{ for all }m\in\mathcal M\bigr\},\\
\nu&=|I|=\prod_{m\in\mathcal M}(m-1).
\end{align*}
Thus $I$ is the Cartesian product of the sets $\{1,\dots,m-1\}$, $m\in\mathcal M$, indexed by the finite set $\mathcal M$; its elements are the families $a=(a_m)_{m\in\mathcal M}$, where we write $a_m=a(m)\in\{1,\dots,m-1\}$ for the coordinate of $a$ indexed by $m$.
(If $\mathcal M=\{m_1,\dots,m_p\}$ with $m_1<\dots<m_p$, then $I$ may be identified with $\{1,\dots,m_1-1\}\times\dots\times\{1,\dots,m_p-1\}$ and $a_{m_j}$ is the $j$-th entry of $a$.)
If $\mathcal M=\emptyset$, then $I$ consists of the empty family only, so $\nu=1$.
Define symmetric $I\times I$ matrices $D_m$ for $m\in\mathcal M\cup\{2,\infty\}$ as follows: $D_2=0$, $D_\infty=2\cdot\mathrm{Id}$, and for $m\in\mathcal M$
\[
(D_m)_{a,b}=\begin{cases}1&\text{if }|a_m-b_m|=1\text{ and }a_{m'}=b_{m'}\text{ for all }m'\in\mathcal M\setminus\{m\},\\0&\text{otherwise.}\end{cases}
\]

\begin{rex}\label{rem:example}
Let $S=\{s_1,s_2,s_3,s_4\}$ with $m_{s_1s_2}=3$, $m_{s_2s_3}=4$, $m_{s_3s_4}=\infty$ and $m_{s_is_j}=2$ for $|i-j|\ge2$, so that all four labels $2,3,4,\infty$ occur.
Then $\mathcal M=\{3,4\}$ and
\[
I=\{1,2\}\times\{1,2,3\}=\{(a_3,a_4)\mid a_3\in\{1,2\},\ a_4\in\{1,2,3\}\},\qquad \nu=6 .
\]
Order $I$ lexicographically: $(1,1),(1,2),(1,3),(2,1),(2,2),(2,3)$.
By definition, $(D_3)_{a,b}=1$ if and only if $a_4=b_4$ and $\{a_3,b_3\}=\{1,2\}$, and $(D_4)_{a,b}=1$ if and only if $a_3=b_3$ and $|a_4-b_4|=1$. With respect to this ordering,
\[
D_3=\begin{pmatrix}
0&0&0&1&0&0\\
0&0&0&0&1&0\\
0&0&0&0&0&1\\
1&0&0&0&0&0\\
0&1&0&0&0&0\\
0&0&1&0&0&0
\end{pmatrix}, \qquad
D_4=\begin{pmatrix}
0&1&0&0&0&0\\
1&0&1&0&0&0\\
0&1&0&0&0&0\\
0&0&0&0&1&0\\
0&0&0&1&0&1\\
0&0&0&0&1&0
\end{pmatrix},
\]
Moreover $D_2=0$ and $D_\infty=2\,\mathrm{Id}_6$.
The idea is now to construct an integral Coxeter system $(\Wt,\St)$ such that $\St$ contains $\nu = 6$ copies of each $s_i$ and the matrices $D_m$ record which of these copies are joined in the Coxeter graph of $\Wt$. We will then show that there is a homomorphism $W \rightarrow \Wt$ which maps each $s_i$ to the product of the copies of $s_i$ in $\Wt$.
\end{rex}

Put $\lambda_m=2\cos(\pi/m)$ for integers $m\ge2$ and $\lambda_\infty=2$, so that $\lambda_{m_{st}}=-2B(\alpha_s,\alpha_t)$ for $s\ne t$.
Define $v\in\R^I$ by
\[
v_a=\prod_{m\in\mathcal M}\sin\Big(\frac{a_m\pi}{m}\Big)\qquad(a\in I).
\]
Since $1\le a_m\le m-1$, we have $v_a>0$ for all $a\in I$.

\begin{lemma}\label{lem:eigen}
$D_mv=\lambda_mv$ for all $m\in\mathcal M\cup\{2,\infty\}$.
\end{lemma}

\begin{proof}
For $m\in\{2,\infty\}$ this is clear. Let $m\in\mathcal M$ and $a\in I$, and put $\theta=a_m\pi/m$ and $v'_a=\prod_{m'\in\mathcal M\setminus\{m\}}\sin(a_{m'}\pi/m')$. The indices $b$ with $(D_m)_{a,b}=1$ are obtained from $a$ by changing $a_m$ into $a_m\pm1$, provided $1\le a_m\pm1\le m-1$. Since $\sin(0)=\sin(\pi)=0$, we obtain in all cases
\[
(D_mv)_a=\big(\sin(\theta-\pi/m)+\sin(\theta+\pi/m)\big)v'_a=2\cos(\pi/m)\sin(\theta)\,v'_a=\lambda_mv_a. \qedhere
\]
\end{proof}

Let $(\Wt,\St)$ be the Coxeter system with $\St=\{r_{s,a}\mid s\in S,\ a\in I\}$ (with $|S|\cdot\nu$ elements) and Coxeter matrix $\widetilde m$ given, for $(s,a)\ne(t,b)$, by
\[
\widetilde m\big(r_{s,a},r_{t,b}\big)=\begin{cases}2&\text{if }s=t,\\ 2,\ 3,\ \infty&\text{if }s\ne t\text{ and }(D_{m_{st}})_{a,b}=0,\ 1,\ 2,\text{ respectively.}\end{cases}
\]
This is well defined, since all $D_m$ are symmetric with entries in $\{0,1,2\}$.
The Coxeter system $(\Wt,\St)$ is integral in the sense of Section~\ref{sec:crys}, and we use the notation introduced there; we write $\at_{s,a}$ for the simple root of $r_{s,a}$ and $\Nt(\tilde w)$ for the set $N(\tilde w)$ of Theorem~\ref{thm:roots}(e) for $(\Wt,\St)$.
By definitions of the form $\Bt$ and the matrices $D_m$,
\begin{equation}\label{eq:Bt}
2\Bt(\at_{s,a},\at_{t,b})=\begin{cases}2\delta_{a,b}&\text{if }s=t,\\-(D_{m_{st}})_{a,b}&\text{if }s\ne t.\end{cases}
\end{equation}
For $s\in S$ put
\[
q_s=\prod_{a\in I}r_{s,a}\in\Wt ;
\]
the factors pairwise commute, so the product does not depend on the order of the factors.

\begin{remark}
If $s \neq t$ and $m_{st} = m$, then in the Coxeter graph of $\Wt$ we have:
\begin{itemize}
    \item \begin{tikzcd}
r_{s,a} \arrow[r, no head,  ""] & r_{t,b}
\end{tikzcd} if $(D_m)_{a,b} = 1$
    \item \begin{tikzcd}
r_{s,a} \arrow[r, no head,  "\infty"] & r_{t,b}
\end{tikzcd} if $(D_m)_{a,b} = 2$
\item there is no edge between
$r_{s,a}$ and $r_{t,b}$ (i.e. they commute) if $(D_m)_{a,b} = 0$
\end{itemize}
\end{remark}

\begin{example}
We proceed with the example from Remark \ref{rem:example}. $\Wt$ has $24$ simple reflections and copies of the same $s_i$ commute. The copies of $s_i$ are given by $r_{s_i, a}$, where $a = (a_3, a_4) \in I$. We exhibit some of the structure of the Coxeter graph of $\Wt$.
\begin{itemize}
    \item We have \begin{tikzcd}
s_1 \arrow[r, no head,  ""] & s_2
\end{tikzcd} in the Coxeter graph of $W$. This yields $6$ copies of type $A_2$ in the Coxeter graph of $\Wt$:
\begin{center}
    \begin{tikzcd}
r_{s_1,(1,j)} \arrow[r, no head,  ""] & r_{s_2,(2,j)}
\end{tikzcd} and
\begin{tikzcd}
r_{s_1,(2,j)} \arrow[r, no head,  ""] & r_{s_2,(1,j)}
\end{tikzcd},
\end{center}
where $j = 1,2,3$.
That is, $D_3$ is the adjacency matrix of the path $A_2$ in the coordinate $a_3$. In each copy of $A_2$, the product of the two generators has order $3$, while products of different copies commute, so $(q_{s_1}q_{s_2})^3 = 1$, matching the relation $m_{s_1 s_2} = 3$ in $W$. 
\item We have \begin{tikzcd}
s_2 \arrow[r, no head,  "4"] & s_3
\end{tikzcd} in the Coxeter graph of $W$. This yields $4$ copies of type $A_3$ in the Coxeter graph of $\Wt$:
\begin{center}
    \begin{tikzcd}
r_{s_2,(i,1)} \arrow[r, no head,  ""] & r_{s_3,(i,2)} \arrow[r, no head,  ""] & r_{s_2,(i,3)}
\end{tikzcd} and
\begin{tikzcd}
r_{s_3,(i,3)} \arrow[r, no head,  ""] & r_{s_2,(i,2)} \arrow[r, no head,  ""] & r_{s_3,(i,1)}
\end{tikzcd},
\end{center}
where $i = 1,2$.
That is, $D_4$ is the adjacency matrix of the path $A_3$ in the coordinate $a_4$. In each copy of $A_3$, the product of the generators is a Coxeter element of type $A_3$, hence of order $4$. So $(q_{s_2}q_{s_3})^4 = 1$, matching the relation $m_{s_2 s_3} = 4$ in $W$. Note that label $4$ itself never appears in $\Wt$.
\end{itemize}
The vector $v$ defined above is
$v=\tfrac{\sqrt3}{2}\bigl(\tfrac{\sqrt2}{2},1,\tfrac{\sqrt2}{2},\tfrac{\sqrt2}{2},1,\tfrac{\sqrt2}{2}\bigr)$, and one checks directly that $D_3v=v=2\cos(\pi/3)\,v$ and $D_4v=\sqrt2\,v=2\cos(\pi/4)\,v$.
\end{example}

\begin{lemma}\label{lem:hom}
There is a unique group homomorphism $F\colon W\to\Wt$ with $F(s)=q_s$ for all $s\in S$.
\end{lemma}

\begin{proof}
By Dyck's theorem and by the Coxeter presentation of $W$ we have to show that $q_s^2=1$ for $s\in S$ and $(q_sq_t)^{m_{st}}=1$ for $s\ne t$ with $m_{st}<\infty$.
Since the $r_{s,a}$ ($a\in I$) are pairwise commuting involutions, $q_s^2=1$.
If $m_{st}=2$, then $D_2=0$ shows that every $r_{s,a}$ commutes with every $r_{t,b}$, so $q_s$ and $q_t$ commute and $(q_sq_t)^2=1$.

Let $m=m_{st}\in\mathcal M$. Let $\Delta$ be the Coxeter graph of the standard parabolic subgroup of $\Wt$ generated by $\{r_{s,a},r_{t,a}\mid a\in I\}$: its vertices are these generators, and two of them are joined if and only if their label is $3$ (no label $\infty$ occurs among them).
Put $I'=\prod_{m'\in\mathcal M\setminus\{m\}}\{1,\dots,m'-1\}$, and for $1\le i\le m-1$ and $a'\in I'$ let $(i,a')\in I$ be the family $b$ with $b_m=i$ and $b_{m'}=a'_{m'}$ for $m'\in\mathcal M\setminus\{m\}$.
By definition of $\widetilde{m}$, two distinct vertices $r_{s,a}, ~r_{s,b}$ are never joined, and neither are $r_{t,a},~r_{t,b}$.
By definition of $D_m$, the vertices $r_{s,(i,a')}$ and $r_{t,(j,b')}$ are joined if and only if $a'=b'$ and $|i-j|=1$. 

For $u \in \{s, t\}$ let $u'$ denote the other element of $\{s,t\}$, and for $a'\in I'$ let $\Delta_{u,a'}$ be the full subgraph on the vertices
\[
r_{u_1,(1,a')},\ r_{u_2,(2,a')},\ \dots,\ r_{u_{m-1},(m-1,a')},\qquad\text{where }u_i=u\text{ for $i$ odd and }u_i=u'\text{ for $i$ even}.
\]
Then $\Delta$ is the disjoint union of the $2|I'|$ subgraphs $\Delta_{u,a'}$, each of which is a path with $m-1$ vertices, i.e.\ a Coxeter graph of type $A_{m-1}$, and there are no edges between different subgraphs.
Consequently, the subgroups $\Wt_{u,a'}$ generated by the vertices of $\Delta_{u,a'}$ pairwise commute, and each of them is a Coxeter group of type $A_{m-1}$ with respect to these vertices \cite[Section~5.5]{Hum90}.
Let $q_{s,u,a'}$ (resp.\ $q_{t,u,a'}$) be the product of the vertices of $\Delta_{u,a'}$ of the form $r_{s,a}$ (resp.\ $r_{t,a}$). Then
\[
q_sq_t=\prod_{u,a'}q_{s,u,a'}\,q_{t,u,a'},
\]
where the factors $q_{s,u,a'}q_{t,u,a'}\in\Wt_{u,a'}$ pairwise commute.
Each $q_{s,u,a'}q_{t,u,a'}$ is a product of all simple generators of the Coxeter group $\Wt_{u,a'}$ of type $A_{m-1}$, each taken once, that is, a Coxeter element of $\Wt_{u,a'}$; it has order $m$, the Coxeter number of type $A_{m-1}$ \cite[Sections~3.16--3.18]{Hum90}.
Hence $(q_sq_t)^m=1$.
If $m_{st}=\infty$ there is nothing to show.
\end{proof}

\subsection{The projection}
Let $\pr\colon\Vt\to V$ be the linear map with
\[
\pr(\at_{s,a})=v_a\,\alpha_s\qquad(s\in S,\ a\in I).
\]
Let $V_+=\{\sum_{s\in S}c_s\alpha_s\in V\setminus\{0\}\mid c_s\ge0\text{ for all }s\in S\}$; note that $V_+\cap(-V_+)=\emptyset$.

\begin{lemma}\label{lem:proj}
The map $\pr$ has the following properties.
\begin{enumerate}[label=\textup{(\alph*)}]
\item $\pr\circ F(w)=w\circ\pr$ for all $w\in W$.
\item Let $\beta\in\Phit$. Then $\beta\in\Phit^+$ if and only if $\pr(\beta)\in V_+$, and $\beta\in\Phit^-$ if and only if $\pr(\beta)\in-V_+$.
\end{enumerate}
\end{lemma}

\begin{proof}
(a) It suffices to treat $w=s\in S$, and to evaluate both sides on the basis vectors $\at_{t,b}$.
Since the $r_{s,a}$ ($a\in I$) are reflections whose roots satisfy $\Bt(\at_{s,a},\at_{s,a'})=\delta_{a,a'}$ by \eqref{eq:Bt}, induction on the number of factors gives
\[
F(s)(\tilde v)=q_s(\tilde v)=\tilde v-\sum_{a\in I}2\Bt(\at_{s,a},\tilde v)\,\at_{s,a}\qquad(\tilde v\in\Vt).
\]
For $\tilde v = \at_{s,b}$ this gives $q_s(\at_{s,b})=\at_{s,b} - 2\at_{s,b} = -\at_{s,b}$, so 
$$(\pr \circ F(s)) (\at_{s,b})=\pr(q_s(\at_{s,b}))=-\pr(\at_{s,b})=-v_b\alpha_s=s(\pr(\at_{s,b})).$$
For $\tilde v = \at_{t,b}$ with $s \neq t$, we obtain
\begin{eqnarray*}
    (\pr \circ F(s)) (\at_{t,b})  =  
    \pr\big(q_s(\at_{t,b})\big)
    & = & 
    \pr\big( \at_{t,b} -\sum_{a \in I} 2\Bt(\at_{s,a}, \at_{t,b}) \at_{s,a}\big)\\
    & \stackrel{\eqref{eq:Bt}}{=} &
    \pr\big( \at_{t,b} +\sum_{a \in I} (D_{m_{st}})_{a,b} \at_{s,a}\big)\\
    & = & 
    v_b \alpha_t + \sum_{a \in I} (D_{m_{st}})_{a,b} v_a \alpha_s\\
    & = & 
    v_b \alpha_t + \sum_{a \in I} (D_{m_{st}})_{b,a} v_a \alpha_s\\
    & = &
    v_b \alpha_t +  (D_{m_{st}} v)_{b}\alpha_s\\
    & \stackrel{\text{(Lemma~\ref{lem:eigen})}}{=} &
    v_b \alpha_t +  (\lambda_{m_{st}} v)_{b}\alpha_s\\
    & = &  v_b\big(\alpha_t+\lambda_{m_{st}}\alpha_s\big)
\end{eqnarray*}
and $s(\pr(\at_{t,b}))=v_b\big(\alpha_t-2B(\alpha_s,\alpha_t)\alpha_s\big)=v_b(\alpha_t+\lambda_{m_{st}}\alpha_s)$.

(b) Let $\beta\in\Phit^+$, say $\beta=\sum_{t,b}c_{t,b}\at_{t,b}$ with all $c_{t,b}\ge0$ and some $c_{t,b}>0$. Then $\pr(\beta)=\sum_t\big(\sum_bc_{t,b}v_b\big)\alpha_t\in V_+$, since all $v_b>0$. If $\beta\in\Phit^-$, then $\pr(\beta)=-\pr(-\beta)\in-V_+$. The assertion now follows from $\Phit=\Phit^+\sqcup\Phit^-$ and $V_+\cap(-V_+)=\emptyset$.
\end{proof}

\subsection{Components of reflections}

\begin{definition}\label{def:comp}
For $g\in T$ let $\Comp(g)=\{\rho\in\Tt\mid \pr(\at_\rho)\in\R_{>0}\,\alpha_g\}$ be the set of \emph{components} of $g$.
\end{definition}

\begin{lemma}\label{lem:comp}
Let $g\in T$ and $w\in W$.
\begin{enumerate}[label=\textup{(\alph*)}]
\item $|\Comp(g)|=\nu$, the elements of $\Comp(g)$ pairwise commute, and their product is $F(g)$.
\item If $g'\in T$ and $g'\ne g$, then $\Comp(g)\cap\Comp(g')=\emptyset$.
\item $F(w)\Comp(g)F(w)^{-1}=\Comp(wgw^{-1})$.
\item If $g\in N(w)$, then $\Comp(g)\subseteq\Nt(F(w))$. If $g\notin N(w)$, then $\Comp(g)\cap\Nt(F(w))=\emptyset$.
\end{enumerate}
\end{lemma}

\begin{proof}
(a) Choose $z\in W$ and $s\in S$ with $g=zsz^{-1}= s_{z(\alpha_s)}$. Then $z(\alpha_s)=\pm\alpha_g$ by Theorem~\ref{thm:roots}(d). Replacing $z$ by $zs$ if necessary, we may assume $z(\alpha_s)=\alpha_g$.
For $a\in I$ put $\beta_a=F(z)(\at_{s,a})\in\Phit$ and $\rho_a=F(z)r_{s,a}F(z)^{-1}\in\Tt$, so that $\rho_a=s_{\beta_a}$.
By Lemma~\ref{lem:proj}(a), $\pr(\beta_a)=z(\pr(\at_{s,a}))=v_a\,z(\alpha_s)=v_a\alpha_g\in V_+$, hence $\beta_a\in\Phit^+$ by Lemma~\ref{lem:proj}(b), so $\beta_a=\at_{\rho_a}$ and $\rho_a\in\Comp(g)$.
Conversely, let $\rho\in\Comp(g)$ and put $\beta'=F(z)^{-1}(\at_\rho)\in\Phit$. By Lemma~\ref{lem:proj}(a), $\pr(\beta')=z^{-1}(\pr(\at_\rho))\in\R_{>0}\,z^{-1}(\alpha_g)=\R_{>0}\,\alpha_s\subseteq V_+$, so $\beta'\in\Phit^+$ by Lemma~\ref{lem:proj}(b).
Write $\beta'=\sum_{t,b}c_{t,b}\at_{t,b}$; by Lemma~\ref{lem:integral} the $c_{t,b}$ are nonnegative integers.
The $\alpha_t$-coordinate of $\pr(\beta')$ is $\sum_bc_{t,b}v_b$, which vanishes for $t\ne s$; since all $v_b>0$, we get $c_{t,b}=0$ for $t\ne s$.
By Lemma~\ref{lem:integral} and \eqref{eq:Bt}, $1=\Bt(\beta',\beta')=\sum_bc_{s,b}^2$, so $\beta'=\at_{s,a}$ for some $a\in I$.
Hence 
$$
\at_\rho=F(z)(\beta') = F(z)(\at_{s,a}) =\beta_a
$$ 
and $\rho=\rho_a$.
Therefore $\Comp(g)=\{\rho_a\mid a\in I\}$.
Since $\Bt(\beta_a,\beta_b)=\Bt(\at_{s,a},\at_{s,b})=\delta_{a,b}$, the roots $\beta_a$ are pairwise distinct, so $|\Comp(g)|=\nu$.
The $\rho_a$ pairwise commute because the $r_{s,a}$ do, and 
$$\prod_a\rho_a=F(z)q_sF(z)^{-1}=F(zsz^{-1})=F(g).$$

(b) If $\rho\in\Comp(g)\cap\Comp(g')$, then $\pr(\at_\rho)$ is a positive multiple of both $\alpha_g$ and $\alpha_{g'}$, so $\alpha_{g'}\in\R_{>0}\,\alpha_g$. Theorem~\ref{thm:roots}(c) gives $\alpha_{g'}=\alpha_g$, whence $g'=g$.

(c) Let $\rho\in\Comp(g)$ and $\beta=\at_\rho$. Write $w(\alpha_g)=\varepsilon\,\alpha_{wgw^{-1}}$ with $\varepsilon\in\{\pm1\}$ (Theorem~\ref{thm:roots}(d)). The element $\varepsilon F(w)(\beta)$ is a root with
\[
\pr(\varepsilon F(w)(\beta))=\varepsilon\, w(\pr(\beta))\in\R_{>0}\,\varepsilon\,w(\alpha_g)=\R_{>0}\,\alpha_{wgw^{-1}}\subseteq V_+,
\]
so it is a positive root by Lemma~\ref{lem:proj}(b). As $F(w)\rho F(w)^{-1}=s_{F(w)(\beta)}=s_{\varepsilon F(w)(\beta)}$, it follows that $F(w)\rho F(w)^{-1}\in\Comp(wgw^{-1})$.
Thus conjugation by $F(w)$ maps $\Comp(g)$ injectively into $\Comp(wgw^{-1})$, and equality holds since both sets have $\nu$ elements by (a).

(d) Let $\rho\in\Comp(g)$. Then $F(w)^{-1}(\at_\rho)\in\Phit$ and, by Lemma~\ref{lem:proj}(a), $\pr\big(F(w)^{-1}(\at_\rho)\big)=w^{-1}(\pr(\at_\rho))\in\R_{>0}\,w^{-1}(\alpha_g)$.
If $g\in N(w)$, then $w^{-1}(\alpha_g)\in\Phi^-$, so this vector lies in $-V_+$, and $F(w)^{-1}(\at_\rho)\in\Phit^-$ by Lemma~\ref{lem:proj}(b), i.e.\ $\rho\in\Nt(F(w))$.
If $g\notin N(w)$, then $w^{-1}(\alpha_g)\in\Phi^+$, and the same argument gives $F(w)^{-1}(\at_\rho)\in\Phit^+$, i.e.\ $\rho\notin\Nt(F(w))$.
\end{proof}

\subsection{Expansions}

\begin{definition}\label{def:expansion}
Let $A=(t_1,\dots,t_k)$ be a tuple of reflections of $W$. An \emph{expansion} of $A$ is a tuple of length $\nu k$ obtained by concatenating, for $i=1,\dots,k$ in this order, the $\nu$ elements of $\Comp(t_i)$ listed in some order.
\end{definition}

By Lemma~\ref{lem:comp}(a), every expansion $\widetilde A$ of $A=(t_1,\dots,t_k)$ consists of reflections of $\Wt$, the product of its entries is $F(t_1)\cdots F(t_k)$, and $\langle\widetilde A\rangle=\big\langle\bigcup_{i=1}^k\Comp(t_i)\big\rangle$.

\begin{lemma}\label{lem:expansion}
Let $c$ be a Coxeter element of $(W,S)$.
\begin{enumerate}[label=\textup{(\alph*)}]
\item $F(c)$ is a Coxeter element of $(\Wt,\St)$, and $\ell_{\Tt}(F(c))=\nu|S|$.
\item Let $x\in[1,c]_T$, put $k=\lT(x)$, and let $A\in\Red_T(x)$. Then every expansion $\widetilde A$ of $A$ lies in $\Red_{\Tt}(F(x))$, and $F(x)\in[1,F(c)]_{\Tt}$. In particular $\ell_{\Tt}(F(x))=\nu k$.
\end{enumerate}
\end{lemma}

\begin{proof}
(a) Write $c=s_1\cdots s_n$ with $S=\{s_1,\dots,s_n\}$. Then $F(c)=q_{s_1}\cdots q_{s_n}$ is a product of all elements $r_{s,a}$ of $\St$, each taken exactly once, so it is a Coxeter element of $(\Wt,\St)$. By Lemma~\ref{thm:lengthc}, $\ell_{\Tt}(F(c))=|\St|=n\nu$.

(b) Let $A=(t_1,\dots,t_k)$. By Lemma~\ref{thm:lengthc} and Lemma~\ref{lem:elementary}(c) there are $t_{k+1},\dots,t_n\in T$ with $(t_1,\dots,t_n)\in\Red_T(c)$.
Let $\widetilde Y$ be the concatenation of $\widetilde A$ with an expansion of $(t_{k+1},\dots,t_n)$. Then $\widetilde Y$ is a tuple of $n\nu$ reflections of $\Wt$ with product $F(t_1)\cdots F(t_n)=F(c)$, so $\widetilde Y\in\Red_{\Tt}(F(c))$ by (a).
The tuple $\widetilde A$ is the prefix of length $\nu k$ of $\widetilde Y$, and the product of its entries is $F(t_1\cdots t_k)=F(x)$. Lemma~\ref{lem:elementary}(b), applied in $(\Wt,\St)$, gives $\widetilde A\in\Red_{\Tt}(F(x))$ and $F(x)\le_{\Tt}F(c)$.
\end{proof}

\subsection{Transfer of lengths}

\begin{proposition}\label{prop:transfer}
Let $R\subseteq T$, let $L=\langle R\rangle$, and let $H=\big\langle\bigcup_{g\in R}\Comp(g)\big\rangle\subseteq\Wt$. Then $F(L)\subseteq H$, and
\[
\nu\cdot\ell_{\chi(L)}(w)\le\ell_{\chi(H)}(F(w))\qquad\text{for all }w\in L .
\]
\end{proposition}

\begin{proof}
For $g\in R$ we have $F(g)\in H$ by Lemma~\ref{lem:comp}(a); hence $F(L)\subseteq H$.
Let $t\in L\cap T$. By Theorem~\ref{thm:dyer}(c) there are $h\in L$ and $g\in R$ with $t=hgh^{-1}$, and Lemma~\ref{lem:comp}(c) gives $\Comp(t)=F(h)\Comp(g)F(h)^{-1}\subseteq H$.

Now let $w\in L$. For every $t\in N(w)\cap L$ we have $t\in L\cap T$, hence $\Comp(t)\subseteq H$, and $\Comp(t)\subseteq\Nt(F(w))$ by Lemma~\ref{lem:comp}(d). Thus
\[
\bigcup_{t\in N(w)\cap L}\Comp(t)\subseteq\Nt(F(w))\cap H .
\]
By Lemma~\ref{lem:comp}(a),(b) the union on the left is disjoint and each $\Comp(t)$ has $\nu$ elements. Hence $\nu\cdot|N(w)\cap L|\le|\Nt(F(w))\cap H|$.
Since $w\in L$ and $F(w)\in H$, Theorem~\ref{thm:dyer}(b), applied to $L\subseteq W$ and to $H\subseteq\Wt$, gives $\ell_{\chi(L)}(w)=|N(w)\cap L|$ and $\ell_{\chi(H)}(F(w))=|\Nt(F(w))\cap H|$.
\end{proof}

\section{Proof of the Main Theorem}\label{sec:proof}

Let $(W,S)$ be a Coxeter system of finite rank, $c\in W$ be a Coxeter element, and $x\in[1,c]_T$. Put $k=\lT(x)$ and let $(t_1,\dots,t_k)\in\Red_T(x)$.
Let $\nu$, $(\Wt,\St)$, $F$ and $\Comp$ be as in Section~\ref{sec:unfolding}.
If $k=0$, then $x=1$, and there is nothing to prove. Assume $k\ge1$.

\begin{proposition}\label{prop:key}
Let $A,A'\in\Red_T(x)$, and let $L$ be the subgroup generated by the entries of $A$ and of $A'$.
Then $|\chi(L)|\le2k$ and $\ell_{\chi(L)}(x)=k$.
\end{proposition}

\begin{proof}
Let $R\subseteq T$ be the set of entries of $A$ and $A'$, so $L=\langle R\rangle$ and $|R|\le2k$; Theorem~\ref{thm:dyer}(c) gives $|\chi(L)|\le2k$.
Let $H=\big\langle\bigcup_{g\in R}\Comp(g)\big\rangle\subseteq\Wt$, and let $\widetilde A$ and $\widetilde{A'}$ be expansions of $A$ and $A'$.
By Lemma~\ref{lem:expansion}, $\widetilde A,\widetilde{A'}\in\Red_{\Tt}(F(x))$ and $F(x)\in[1,F(c)]_{\Tt}$, where $F(c)$ is a Coxeter element of the integral Coxeter system $(\Wt,\St)$ of finite rank.
By Theorem~\ref{thm:hk}, $\widetilde{A'}$ lies in the Hurwitz orbit of $\widetilde A$, so $\langle\widetilde A\rangle=\langle\widetilde{A'}\rangle$ by Lemma~\ref{lem:elementary}(a).
Since $\langle\widetilde A\rangle$ and $\langle\widetilde{A'}\rangle$ are generated by the components of the entries of $A$ and of $A'$, respectively, we get
\[
H=\big\langle\langle\widetilde A\rangle\cup\langle\widetilde{A'}\rangle\big\rangle=\langle\widetilde A\rangle .
\]
Proposition~\ref{prop:crys}, applied to $y=F(x)$, $q=\nu k$ and $(g_1,\dots,g_q)=\widetilde A$, yields $\ell_{\chi(H)}(F(x))=\nu k$.
Since $x\in L$, Proposition~\ref{prop:transfer} gives $\nu\cdot\ell_{\chi(L)}(x)\le\ell_{\chi(H)}(F(x))=\nu k$, that is, $\ell_{\chi(L)}(x)\le k$.
On the other hand $\ell_{\chi(L)}(x)\ge\lT(x)=k$ by \eqref{eq:lengthineq}.
\end{proof}

\begin{proof}[Proof of the Main Theorem]
\emph{(iii)} Let $A,A'\in\Red_T(x)$ and let $L$ be the subgroup generated by the entries of $A$ and $A'$.
By Proposition~\ref{prop:key}, $\chi(L)$ is finite and $\ell_{\chi(L)}(x)=\lT(x)=k$.
By Lemma~\ref{lem:parabolic}(a) there are pairwise distinct $r_1,\dots,r_k\in\chi(L)$ with $x=r_1\cdots r_k$, and by Lemma~\ref{lem:parabolic}(b) both $A$ and $A'$ lie in the Hurwitz orbit of $(r_1,\dots,r_k)$.
Hence $A$ and $A'$ lie in the same Hurwitz orbit.

\emph{(i)} Put $W_0=\langle t_1,\dots,t_k\rangle$. Let $t\in T$ with $t\leT x$. By Lemma~\ref{lem:elementary}(d) there is $A'\in\Red_T(x)$ whose first entry is $t$. By (iii), $A'$ lies in the Hurwitz orbit of $(t_1,\dots,t_k)$, so $t\in\langle A'\rangle=W_0$ by Lemma~\ref{lem:elementary}(a). Hence $P(x)\subseteq W_0$. Conversely, $t_i\leT x$ for $1\le i\le k$ by Lemma~\ref{lem:elementary}(d), so $W_0\subseteq P(x)$.

\emph{(ii)} Apply Proposition~\ref{prop:key} with $A=A'=(t_1,\dots,t_k)$. Then $L=W_0=P(x)$ by (i), and $\ell_{\chi(P(x))}(x)=k$.
By Lemma~\ref{lem:parabolic}(a) there are pairwise distinct $r_1,\dots,r_k\in\chi(P(x))$ with $x=r_1\cdots r_k$.
Since $P(x)=\langle t_1,\dots,t_k\rangle$ is generated by $k$ reflections, Theorem~\ref{thm:dyer}(c) gives $|\chi(P(x))|\le k$.
Hence $\chi(P(x))=\{r_1,\dots,r_k\}$, so $|\chi(P(x))|=k$ and $x$ is a Coxeter element of $(P(x),\chi(P(x)))$.
\end{proof}

\bibliography{prefixes}
\bibliographystyle{amsplain}

\end{document}